\documentclass[11pt, twoside]{amsart}

\usepackage{amssymb, amsmath, amsthm, mathrsfs,}
\usepackage{graphicx, epsfig, color,wrapfig}
\usepackage{tikz}
\usepackage{hyperref}
\usepackage[mathscr]{euscript}
\usepackage{enumerate}
\usepackage{fullpage}
\graphicspath{{figures/}} 

\usepackage{cancel} 
\usepackage{hyperref} 

\hypersetup{
    colorlinks,%
    citecolor=black,%
    filecolor=black,%
    linkcolor=black,%
    urlcolor=blue
}

\usepackage{todonotes}

\usepackage{lipsum}
\usepackage{floatrow}
\usepackage[latin1]{inputenc}
\usepackage[all,cmtip]{xy}

\usepackage{chngpage}

\usepackage{booktabs}

\newtheorem{theorem}{Theorem}[section]
\newtheorem{definition}[theorem]{Definition}

\newtheorem{proposition}[theorem]{Proposition}

\newtheorem{lemma}[theorem]{Lemma}
\newtheorem{corollary}[theorem]{Corollary}

\newtheorem*{theorem*}{Theorem}
\newtheorem*{lemma*}{Lemma}
\newtheorem*{exercise*}{Exercise}

\theoremstyle{definition} 
\newtheorem*{definition*}{Definition}
\newtheorem*{notation*}{Notation}
\newtheorem*{example*}{Example}
\newtheorem*{question*}{Question} 

\title[Ranks of Maps of Vector Bundles]{Ranks of Maps of Vector Bundles}

\author[M. Teixidor]{Montserrat Teixidor {\textrm i} Bigas} \address{Department of Mathematics,  Tufts University, 177 College Ave,  Medford, MA 02155, U.S.A.}

 \email{mteixido@tufts.edu}

\dedicatory{Dedicated to  Peter Newstead for his extraordinary mentorship role.}

\begin{document}
\maketitle

\begin{abstract} We generalize to vector bundles the techniques introduced for line bundles in \cite{LOTZ1}. 
We then use this method to prove the injectivity of the Petri map for vector bundles and the surjectivity of a map related to deformation theory of Poincar\'e sheaves.
\end{abstract}

Many questions about the geometry of curves have a translation in terms of the rank of a map among spaces of sections of line bundles.
The best known may be  the Petri map :
\[ H^0(C,L)\otimes H^0(C, K\otimes L^{-1})\to H^0(C, K).\]
The injectivity of this map determines the non-singularity of the locus of special divisors on a curve.
Similarly, the maximal rank conjecture states that for a general curve and a general line bundle, the natural map of sections 
\[ S^m(H^0(C,L))\to H^0(C, L^m)\]
has maximal rank allowing one to count how many linearly independent hypersurfaces in the ambient space of a given degree contain a curve.

 These conjectures have been proved by  degenerating to a special curve and 
 then showing that the kernel of the map gives rise to a small enough subspace in the limit and therefore it needs to be small too  on the generic curve.
A different approach is to look at the image of the map and  show it is large enough.
This amounts    to proving the linear independence  of a collection of sections  obtained as product of sections of  linear series in the domain.
This method was used  for line bundles in search of  proofs of the maximal rank conjecture in \cite{LOTZ1}, \cite{LOTZ2}.
This is also essentially the strategy in  the tropical setting (see   \cite{JP}, \cite{FJP}).

In this paper, we   generalize to bundles of degree higher than one  techniques that so far have only been used for line bundles (see  section \ref{secGener}).
We then show how they can be used  to solve a couple of questions about vector bundles.
In particular, we provide a new and simpler proof of the injectivity of the Petri map for vector bundles 
satisfying only mild  numerical conditions on  rank, degree and number of sections(see section \ref{secPetri}). 
Note that  injectivity of the Petri map  guarantees existence of a component of the right dimension of the Brill Noether locus for vector bundles on a generic curve.
In particular, this method could be used to fill gaps on the known Brill-Noether geography that keeps track of the values $r,d, k$ for which the Brill -Noether locus $W^k_{r,d}$ is non-empty.

We also look (see section \ref{secPoin})  at the surjectivity of the cup product map between the sections of the canonical bundles of the curve 
and the sections of  the traceless endomorphisms tensored with the canonical..
This surjectivity has important consequences for the  deformation theory  of Poincar\'e bundles (see \cite{BBN}).

It is worth pointing out that the applications in this paper are simple because we only try to show results for a generic vector bundle of a given rank and degree
which in practice means checking the result for one vector bundle of that rank and degree.
Checking for every vector bundle  would lead to  substantially more involved combinatorial issues.

While many proofs in the setting of Brill-Noether Theory for line bundles can be run in parallel using either techniques of limit linear series or tropical curves, 
 no satisfactory theory for higher rank has so far been developed in the tropical setting.
 It would be nice to be able to use tropical techniques to solve the many open problems related to higher rank questions on curves but this does not seem possible at this time.
 
 Acknowledgments:   I would like to thank Peter Newstead for suggesting the problem dealt with in  section \ref{secPoin} and to a careful referee for suggesting clarifications.


\section{Vector bundles on reducible curves and limit linear series.}\label{secGener}

The goal of this paper is to show that on a generic curve of fixed genus, certain maps among spaces of sections of vector bundles have at least a certain rank.
We will use degeneration methods. 
We review here the tools and techniques that we will need to deal with linear series on reducible, semistable,  nodal curves.

Most of the definitions and results we recall in this section make sense for an arbitrary nodal curve of compact type 
(that is, one whose Jacobian is compact or equivalently, whose dual graph has no loops).
We will only use them with  a particular type of reducible curve, namely  generic chains of elliptic curves of a fixed genus $g$,
 introduced by Welters in \cite{W} and  that we have been using systematically  in the context of vector bundles of rank one and higher (\cite{SELB}, \cite{langred}, \cite{CT}).

\begin{definition}[generic chain of elliptic curves]\label{chain} A generic chain of elliptic curves is built as follows:
Let  $C_1,\dots,C_M$ be curves that are either rational or elliptic, with precisely $g$ of them being elliptic,
$P_i, Q_i$ generic points on $C_i$, for $i=1,\dots,M$. Glue $Q_i$ to $P_{i+1},i = 1,\dots,M-1$ to form a nodal curve of genus $g$. 
\end{definition}

Consider a family of curves in which the generic curve is non-singular and the special curve is a chain of elliptic curves.
Given a vector bundle together with a space of sections on the generic fiber, up to some base change and normalizations, we can complete it to a vector bundle on the family.
If we start with a chain of elliptic curves as the central fiber, the new central fiber will still be a chain of elliptic curves.
Tensoring the family with line bundles with support on the components of the central fiber, the vector bundle on the generic curve does not change but on the special curve it does.
In particular, the limit on the special curve of a vector bundle in a family is not unique.
In the case of rank one,  this process allows to obtain an arbitrary distribution of degrees among the components of the reducible limit curve.
This is the basis for the theory of limit linear series introduced by Eisenbud and Harris  \cite{EH1}.
These authors took advantage of the possibility of concentrating all of the degree in one of the irreducible components of the central fiber.
Then for every family $g^{k-1}_d$ (we use $k$ rather than $r+1$ to avoid confusions with the rank)
of linear series and for each components $C_i$ of the reducible limit  curve, they obtained a line bundle $L_i$ of degree $d$ on $C_i$. 
When concentrating all the degree in the component $C_i$, the limit of  the $k$-dimensional  space of sections of the family is also  concentrated on $C_i$.
Therefore, there are distinguished  $k$-dimensional  spaces of sections $V_i$ of $L_i$ for each component $C_i$ of the central fiber.
 Compatibility among the data on the different components give rise to vanishing conditions of the sections of the linear series at the nodes:
\begin{definition}\label{llsrk1} [Limit Linear Series, rank one]\cite{EH1}, 
Given a curve of compact type, a crude  limit linear series of degree $d$ and dimension $k$ is the data of a $g^{k-1}_d$ on each irreducible component
such that if   $Y_1,Y_2$ are irreducible components meeting at a point $P$ and $a_1^i,\dots ,a_k^i$ are the orders of vanishing of the sections of these series at $P$ on $Y_i, i=1,2$,
then $a_j^1+a_{k-j+1}^2\ge d$.
The series is refined if these inequalities are equalities.
\end{definition}
 
A similar approach can be taken for vector bundles, but we need to account for three differences
\begin{itemize}
\item For a line  bundle and an $k$-dimensional space of its sections, the orders of vanishing of the sections at any given point are  $k$ distinct integers $a_1<\dots<a_k$.
For a vector bundle $E$ of rank $r$ and an $k$-dimensional space of its sections $V$, the orders of vanishing of the sections at a  point  do not normally take $k$ distinct  values.
But we can find $k$ integers $0\le a_1\le \dots\le a_k$ such that  a particular value $\alpha$ appears $t$ times in the sequence with $t=\dim (V(-\alpha P)/ V(-(\alpha +1)P))$.
\item If $E$ is a vector bundle of rank $r$ on an irreducible component $C_i$ and $P\in C_i$, then $\deg(E(-P))=\deg(E)-r$.
Modifying a vector bundle by tensoring with a line bundle with support on the central fiber changes the degrees of the restriction of the bundle to the components of the central fiber by multiples of $r$.
We will not be able to assume that the limit on the central fiber of a vector bundle on the generic curve has degree 0 on all components of the central fiber except for one.
But we can still assume that in the limit, all sections are concentrated in one component by taking sufficiently small degree (perhaps negative) on the remaining components.
\item While a line bundle on a curve of compact type is completely determined by its restrictions to the components of the curve, vector bundles depend also on the gluing at the nodes.
\end{itemize}
Taking these issues into account, we consider  the concept of Limit Linear Series for Vector Bundles  (see for example  \cite{Clay}, \cite{duke}).

\begin{definition}\label{lls} [Limit Linear Series, arbitrary rank]
A limit linear series of rank $r$, degree $d$ and dimension $k$ on a chain of ellliptic curves consists of the following
\begin{enumerate}[(a)]
\item A vector bundle $E_i$  of rank $r$ and degree $d_i$ on each component $C_i$ and a $k$-dimensional space of sections $V_i$ of $E_i$.
\item An isomorphism between the projectivization of the fibers $(E_i)_{Q_i}$ and $(E_{i+1})_{P_{i+1}}, i=1,\dots, M-1$;
\item Bases $s^{t}_{Q_i}$, $s^{t}_{P_{i+1}}$, $t=1,\dots,k$ of the vector spaces $V_i$ and $V_{i+1}$.
\item A positive integer $a$.
\end{enumerate}
These data satisfy the following conditions:
\begin{enumerate}[(1)]
\item $\sum _{i=1}^M d_i-r(M-1)a=d$.
\item The orders of vanishing  at $Q_{i}$ and $P_{i+1}$ of the sections of the chosen bases satisfy $ord_{P_{i+1}}s^t_{P_{i+1}}+ord_{Q_i}s_{Q_i}^t\geq a$ for all $t$.
\item Sections of the vector bundles $E_i(-aP_i)$ and $E_i(-aQ_i)$ are completely determined by their value at $P_i$ and $Q_i$, respectively.
\end{enumerate}
\end{definition}

This definition of Limit Linear Series amounts to,  at each step, concentrating all of the degree and therefore all of the sections on one component of the central fiber.
It has the advantage that one only needs to deal with one irreducible curve at a time but the disadvantage that  all the sections must be dealt with simultaneously.
This becomes  difficult, if not impossible when trying to show their independence.
To prove our results, we will need  to spread out the degree among the different components. 
The theoretical framework for this alternative point of view was considered in rank one by Osserman \cite{O} and can also be used for vector bundles \cite{tohoku}.
In fact, the data of a limit linear series as defined in \ref{lls} can be used to produce vector bundles and the corresponding spaces of sections 
with degree and sections spread out among the components of the reducible curve.
As the distribution of degrees comes from tensoring  the original vector bundle with line bundles with support on the fibers,
the process changes the degree in  multiples of $r$, so we need to keep the same remainder modulo $r$ for the degree on specific components.
In fact, the restriction to a component will be modified by tensoring with a line bundle with support at the nodes.
The sections from the limit linear series  will give rise to sections of the modified vector bundle if they vanish at the nodes with the right order.
This question was studied in \cite{tohoku}, \cite{documenta}.
More specifically, we have:

\begin{lemma}\label{serdistrdeg}
Given a chain of elliptic curves, a limit linear series  corresponding to the data of vector bundles $E_i$  spaces of sections $V_i$ of $E_i$
and isomorphism of  the projectivized fibers $(E_i)_{Q_i}\cong (E_{i+1})_{P_{i+1}}$
as in \ref{lls} and a distribution of degrees $d'_i$ with $\sum d'_i=d$ and $d'_i$ congruent with $d_i$ modulo $r$,
there exists a vector bundle on the whole curve whose restriction to the component $C_i$ has degree $d'_i$ and a vector space of sections of this vector bundle of dimension $k$.
\end{lemma}
\begin{proof} Given a limit linear series on a chain of elliptic curves, condition (3) means that $d_i=ar+\bar d_i$ where $0\le \bar d_i<r$.
Then, 
\[ d=\sum _{i=1}^M d_i-r(M-1)a=\sum _{i=1}^M \bar d_i+rMa-r(M-1)a=ra+\sum _{i=1}^M \bar d_i \]

As  we assume $d'_i$ congruent with $d_i$ modulo $r$, we can write  $d'_i=a_ir+\bar d_i$.
With this notation, the condition  $\sum d'_i=d$ can be written as  $\sum \bar d_i+(\sum a_i)r=d$.
Define now the vector bundle $E'_i$ and a space of its sections as follows 
\[ E'_i=E_i(-(\sum_{j<i}a_j)rP_i-(\sum_{j>i}a_j)rQ_i),\  V'_i=\{ s\in V_i| \ ord_{P_i}s\ge \sum_{j<i}a_j, \  ord_{Q_i}s\ge \sum_{j>i}a_j\}. \]
Take the  gluing at the nodes induced from those given for the original linear series.

We now check that the degree of $E'_i$ is $d'_i$ as needed: From the definition of $E'_i$, 
\[ \deg (E'_i)=d_i-(\sum_{j\not= i}a_j)r=\bar d_i+ar-(\sum_{j\not= i}a_j)r =d- \sum_{j\not= i}\bar d_j -(\sum_{j\not= i}a_j)r =d-\sum_{j\not= i}d'_j =d'_i \]
as needed.

As the sections of the $V_i$ glued with each other, so do the sections of the newly defined $V'_i$, giving rise to sections of the new vector bundle.
\end{proof}

 In each proof, our strategy to show independence of sections of a product will be to fix the distribution of degrees among the components.
 Then, from a limit linear series, we produce the  sections of the new vector bundle with the prescribed degrees according to Lemma \ref{serdistrdeg}
 and show that these are independent.
If the distribution of degrees has been  chosen judiciously, one can show that for a few or all the sections that are non-zero on the component,
 the coefficients in front of these sections in the linear combination are 0.

We need to build by hand  vector bundles of rank $r$ and degree $d$ with certain number of sections.
For reducible curves $X$, there exist  moduli spaces of vector bundles of given rank and degree that are stable for a given polarization.
These moduli spaces of vector bundles are themselves reducible with components corresponding to the distribution of degrees modulo $r$  among the components of $X$.
Once the distribution of degrees has been fixed, the choice of a semistable bundle on each component of $X$ gives rise to a semistable vector bundle on $X$ 
which is stable so long as possible destabilizing subbundles do not glue with each other. For details, please see \cite{modred}, \cite{modred2}
In particular, we have the following
\begin{theorem} \label{modvbrc}Let  $X=C_1\cup\dots\cup C_g$ be a chain of elliptic curves as in Definition \ref{chain},  $r,d\in {\mathbb Z} $.
Given semistable vector bundles on each components one of which is stable, gluing with arbitrary gluing at the nodes, one obtains a vector bundle on the chain stable by a suitable polarization.
If all the vector bundles on the individual elliptic curves are strictly semistable, the resulting vector bundle is stable provided the destabilizing subbundles do not glue with each other.

The components of the moduli space of vector bundles of rank $r$ and degree $d$ on $X$ that are stable for a given polarization correspond with all the possible distributions of the degrees $d_i$ to $X_i$ modulo $r$.
For a choice of such a component $M$ of the moduli space,  and a component $C_i$ of $X$, if $h_i$ is the greatest common divisor of $r$ and the corresponding $d_i$, 
then the restriction to $C_i$ of a generic element $E\in M$ is a direct sum of $h_i$ indecomposable vector bundles of rank $\frac r{h_i}$ and degree $\frac {d_i}{h_i}$.
 \end{theorem}

In particular, if the degree assigned to a given elliptic component  of our curve is divisible by the rank, we need to consider direct sums of $r$ line bundles. 
So, let us start by looking at a single elliptic curve and the sections of a line bundle on it.

\begin{lemma} \label{lemseclselc}Let $C$ be an elliptic curve, $P,Q\in C$ so that $P-Q$ is not a torsion element (in the group structure of $C$).
Let $L$ be a line bundle of degree $d$ on $C$.
For any $k, 0\le k\le d-1$, there exists up to a constant, a unique section $s_k$ of $L$ such that $s_k$ vanishes at $P$ with order at least $k$ and at $Q$ with order at least $d-k-1$.
The inequalities are in fact equalities   unless $L={\mathcal O}(kP+(d-k)Q)$ when $s_{k-1}=s_k$ or  $L={\mathcal O}((k+1)P+(d-k-1)Q)$ when $s_k=s_{k+1}.$
\end{lemma}
\begin{proof} As $\deg(L(-kP-(d-k-1)Q))=1$, from Rieman-Roch's Theorem, $h^0(C, L(-kP-(d-k-1)Q))=1$.
Choose $s_k\in H^0(C, L(-kP-(d-k-1)Q)), s_k\not=0$.
If $ L\not={\mathcal O}(kP+(d-k)Q), L\not={\mathcal O}((k+1)P+(d-k-1)Q)$, then $h^0(C, L(-kP-(d-k)Q))=0, h^0(C, L(-(k+1)P-(d-k-1)Q))=0$.
Therefore, unless we are in one of these two situations, $s_k$ vanishes at $P$ with order precisely $k$ and at $Q$ with order precisely $d-k-1$.

If there were sections $s_{k_1}, s_{k_2}$ with orders of vanishing at $P, Q$ adding to $d$, then  $L={\mathcal O}(k_iP+(d-k_i)Q), i=1,2$. But
\[ {\mathcal O}(k_1P+(d-k_1)Q)={\mathcal O}(k_2P+(d-k_2)Q)\Rightarrow (k_2-k_1)(P-Q)=0\]
This contradicts our assumption that $P-Q$ is not a torsion point of the curve.
\end{proof}

\begin{corollary} \label{corsecllsec} Let $C$ be an elliptic curve, $P,Q\in C$ so that $P-Q$ is not a torsion element,  $L$  a line bundle of degree $d$ on $C$.
\begin{enumerate}[(a)]
\item If $L={\mathcal O}(aP+(d-a)Q), 0\le u\le a\le u+t\le d-1$, there exists a uniquely determined $t$-dimensional space of sections of $L$ whose orders of vanishing at $P,Q$ are respectively
\[ u, u+1,\dots, a-3, a-2, a, a+1, a+2,\dots,u+t;\ \ d-u-1, d-u-2, \dots, d-a+2, d-a+1, d-a, d-a-2, d-a-3,\dots, d-u-t-1\]
\item If $L\not={\mathcal O}(aP+(d-a)Q)$ for any $a$ and $ 0\le u\le u+t\le d$, there exists a uniquely determined $t$-dimensional space of sections of $L$ whose orders of vanishing at $P,Q$ are respectively
\[ u, u+1,\dots,u+t-1;\ \ d-u-1, d-u-2, \dots, d-u-t\]
\end{enumerate}
\end{corollary}
\begin{proof} The sections $s_k$ defined in Lemma \ref{lemseclselc} have different orders of vanishing at $P$, therefore, they are linearly independent.
If we take $t$ of them, they span a $t$-dimensional space of sections of $L$.
\begin{enumerate}[(a)]
\item If $L={\mathcal O}(aP+(d-a)Q)$, take the vector space spanned by $s_u, s_{u+1},\dots, s_{a-2}, s_{a-1}=s_a, s_{a+1},\dots, s_{u+t}$.
The vanishing conditions at the nodes are satisfied.

Conversely, assume that a space $V$ of sections of dimension $t$ of $L$ has the given vanishing at the nodes.
For $u\le x\le a-1$, the subspace of  $V$ of sections vanishing to order at least $x$ at $P$ has dimension at least $t-x+u$.
while the subspace of those  vanishing to order at least $d-x-1$ at $Q$ has dimension at least $x-u+1$.
As, by assumption,  both subspaces live in a $t$-dimensional space, the two subspaces intersect.
From Lemma \ref{lemseclselc},  $s_x \in V$.

Similarly,  for $a+1\le x\le u+t$, the subspace of sections vanishing to order at least $x$ at $P$ has dimension at least $t-x+u+1$
while the subspace of sections vanishing to order at least $d-x-1$ at $Q$ has dimension at least $x-u$, so again $s_x \in V$.
As the $s_x$ are linearly independent, $V$ is  the span of  $s_u, s_{u+1},\dots, s_{a-2}, s_{a-1}=s_a, s_{a+1},\dots, s_{u+t}$.
\item As in the proof of (a), the sections  $s_u,\dots ,s_{u+t-1}$ span the required subspace of $H^0(C,L)$. 
\end{enumerate}
\end{proof}

\begin{lemma} \label{lemsecvbec} Let $r,d\in {\mathbb Z} , r>0, d\ge 0, d=rd_1+d_2, h=gcd(r,d)$.
Let    $C$  be an elliptic curve, $P\in C$,  $E$  a  vector bundle on $C$ that can be written as direct sum of $h$ indecomposable vector bundles of degree $\frac dh$ and rank  $\frac rh$.
 Then,  there exists a uniquely determined $rk_1+d_2$-dimensional space of sections of $E$ whose orders of vanishing at $P$ are 
\[ \overbrace{d_1, \cdots d_1}^{d_2\rm\ times},  \overbrace{d_1-1 \cdots d_1-1}^{r\rm\ times},  \overbrace{d_1-2 \cdots d_1-2}^{r\rm\ times}, \dots  ,  \overbrace{d_1-k_1 \cdots d_1-k_1}^{r\rm\ times}\]
\end{lemma}
\begin{proof} From Riemann-Roch's Theorem, $h^0(C, E(-\alpha P))=\deg( E(-\alpha P))=r(d_1-\alpha)+d_2$.
\end{proof}

In our applications, we will need to use the limit of the canonical linear series on a curve. 
For ease of reference, we recall here what it looks like on our chains:

\begin{proposition}\label{can}The canonical limit linear series on $C_0$ has line bundles on $C_i$ equal to 
\[ L_i={\mathcal O}(2(i-1)P_i+2(g-i)Q_i) \]
The space of sections on $C_i$ is 
 \[ H^0(L_i(-(i-2)P_i-(2g-2i)Q_i))\oplus H^0(L_i(-(2i-1)P_i-(g-i-1)Q_i)). \]
The unique section whose order of vanishing at $P_i$ and $Q_i$ is $2g-2$  vanishes with order $2(i-1)$ at $P_i$ and $2g-2i$ at $Q_i$.
  \end{proposition}
\begin{proof}  See \cite{EH} Th2.2 for a proof in a more general situation.

The result also follows from the one to one correspondence between limit linear series of degree $d$ and projective dimension $r$ on a general chain of elliptic curves 
and  filling of Young Tableaux of dimensions $(r+1)\times (g-d+r)$ with numbers among the $1,\dots, g$ so that they are strictly increasing on each row and column(see  \cite{documenta}). 
  This description shows that in the case of $W_{2g-2}^{g-1}$ there is a unique limit linear series  and that the restriction of the  line bundle to $C_i$ is of the form $\mathcal{O}(2(i-1)P_i+2(g-i)Q_i)$
  and the space of sections on the component $C_i$ has vanishing at $P_i, Q_i$ respectively as 
  \[  \begin{matrix} i-2 & i-1 &\dots &2i-5&2i-4&2i-2&2i-1&\dots&g+i-2\\               2g-i-1&2g-i-2&\dots &2g-2i+2&2g-2i+1&2g-2i&2g-2i-2&\dots &g-i-1\\  \end{matrix}\]
\end{proof}


\section{Injectivity of the Petri map in higher rank}\label{secPetri}

The Petri map  for line bundles  that we mentioned before  $H^0(C,L)\otimes H^0(C, K\otimes L^{-1})\to H^0(C, K)$, 
controls the geometry of the Brill Noether locus. 
The injectivity of the Petri map guarantees that $L$ gives rise to a non-singular point of the Brill Noether locus of $W_d^r$ for $d=\deg L, r=h^0(C,L)-1$.

A similar Petri map can be defined in higher rank 
\[ P: H^0(C,E)\otimes H^0(C, K\otimes E^*)\to H^0(C, K\otimes E\otimes E^*)\]

While for a generic curve $C$ and every line bundle $L$ on $C$, the Petri map is injective, it is no longer the case that even for the generic curve and every vector bundle, the higher rank Petri map is injective.
On the other hand, the injectivity of the map for a particular vector bundles $E$ guarantees the existence of a component of the expected dimension of the Brill-Noether locus .
Injectivity was proved in a fairly reduced number of cases (small slope or small number of sections) in \cite{BP}, \cite{BGMMN1},\cite{BGMMN2} 
and for rank two with additional conditions on degree, genus and number of sections in \cite{CF}.
A proof for the generic curve with fewer restrictive numerical conditions is provided in \cite{CLT}.
Here we give a different proof using the techniques of section \ref{secGener}.

We prove the following:

\begin{theorem} Fix $r,d,k$.
 Write $d=rd_1+d_2, k=rk_1+k_2, 0\le d_2<r, 0\le k_2<r$.
 Assume that one of the following holds
 \begin{itemize}
 \item  $d_2\ge k_2, \ d_2\not= 0,\ (k_1+1)(g+k_1-d_1-1)\le g-1$
 \item  $d_2= k_2=0,\ k_1(g+k_1-d_1-1)\le g-2$
\item   $d_2< k_2,  (k_1+1)(g+k_1-d_1)\le g-1$
\end{itemize}
Then the generic curve has a component of expected dimension of the Brill-Noether locus of vector bundles of rank $r$ degree $d$ with $k$ sections for which the Petri map is generically injective.
\end{theorem}
\begin{proof} The proof is similar in all cases but to simplify notations, we will assume 
\[ d_2=k_2\not=0,\ (k_1+1)(g+k_1-d_1-1)\le g-1.\]
We construct a  limit linear series on a generic chain of elliptic curves.
\begin{enumerate}[(a)]
 \item  For $ i\le  (k_1+1)(g+k_1-d_1-1)$, write $i=(k_1+1)j_1+j_2,\ 1\le j_2\le k_1+1$ (and therefore $0\le j_1\le g+k_1-d_1-2$).
 If  $j_2\not= k_1+1$, take the vector bundle on $C_i$ to be the direct sum 
 \[ {\mathcal O}_{C_i}((i-j_1+j_2-2)P_i+(d_1-i+j_1-j_2+2)Q_i)^{\oplus r}.\]
  On the first $k_2$ line bundles  ${\mathcal O}_{C_i}((i-j_1+j_2-2)P_i+(d_1-i+j_1-j_2+2)Q_i)$, from Corollary \ref{corsecllsec} (a), we can
  take a space of sections with vanishing orders at $P_i, Q_i$ respectively given as   
 \[  \begin{matrix} i-j_1-2& i-j_1-1&\dots & i-j_1+j_2-4&i-j_1+j_2-2& i-j_1+j_2-1&\dots& i-j_1+k_1-1\\
 d_1-i+j_1+1& d_1-i+j_1&\dots& d _1-i+j_1-j_2+3& d _1-i+j_1-j_2+2&d_1-i+j_1-j_2&\dots& d_1-i+j_1-k_1\end{matrix} \]
 On the remaining $r-k_2$ line bundles, remove the last section.
 \item If $i=(k_1+1)j_1+j_2\le  (k_1+1)(g+k_1-d_1-1), j_2= k_1+1$, take the direct sum
 \[{\mathcal O}_{C_i}((i+j_2-j_1-2)P_i+(d_1+j_1-j_2-i+2)Q_i)^{\oplus k_2}\oplus L_1\oplus \dots\oplus L_{r-k_2}\]
  where the $L_j$ are generic line bundles of degree $d_1$. 
 Take the space of sections as in case (a).
 Note that when $j_2=k_1+1$ the last section is the one whose sum of vanishing at $P_i, Q_i$ is $d_1$.
  This section is omitted when the line bundle corresponds to an index larger than $k_2$.
 From  Corollary \ref{corsecllsec} (b),  this construction is possible.
\item  For $  (k_1+1)(g+k_1-d_1-1)<i\le g-1$, take the vector bundle on $C_i$ to be $ L_1\oplus \dots\oplus L_{r}$ where the $L_j$ are generic line bundles of degree $d_1$. 
 Write $\alpha =g+k_1-d_1-1$. On the first $k_2$ line bundles, take a  space of sections with vanishing orders at $P_i, Q_i$ respectively given as 
 \[  \begin{matrix} i-\alpha-1& i-\alpha&\dots & i-\alpha+k_1-1\\
 d_1-i+\alpha& d_1-i+\alpha-1&\dots&  d_1-i+\alpha -k_1\end{matrix} \]
 on the remaining $r-k_2$ line bundles, remove the last section.
 Again, from  Corollary \ref{corsecllsec} (b),  this construction is possible.
 \item On $C_g$ take $E$ to be the direct sum of $h$ generic stable vector bundles of rank $\frac rh$ and degree  $\frac dh$  where $h$ is the greatest common divisor of $d, r$.
 Take the  space of sections $H^0(C_g, E(-(g-\alpha-1)P_g))$ that is, the space of sections whose orders of vanishing at $P_g$ are 
\[ \overbrace{g-\alpha-1, \cdots g-\alpha-1}^{r\rm\ times},  \overbrace{g-\alpha, \cdots , g-\alpha}^{r\rm\ times}, \dots,   \overbrace{d_1-1, \cdots ,d_1-1}^{r\rm\ times} , 
  \overbrace{d_1, \cdots , d_1}^{d_2\rm\ times} \]
 Using that  $\alpha =g+k_1-d_1-1$, these vanishing can be written as 
 \[ \overbrace{d_1-k_1, \cdots ,d_1-k_1}^{r\rm\ times},  \overbrace{d_1-k_1+1, \cdots ,d_1-k_1+1}^{r\rm\ times},\dots,   \overbrace{d_1-1 , \cdots , d_1-1}^{r\rm\ times} , 
  \overbrace{d_1, \cdots ,d_1}^{d_2\rm\ times} \]
  From Lemma \ref{lemsecvbec}, this is possible. 
\end{enumerate}
The gluing is such that  for $ 2\le i\le  (k_1+1)(g+k_1-d_1-1)$, the line bundles in the decomposition glue with the corresponding line bundles in the decomposition on the previous curve.
After this, the gluing is generic. At the last node, the subspace of sections of $E_g$ of dimension $d_2=k_2$ with maximum vanishing at $P_g$ is glued with the subspace of sections chosen in 
$E_{g-1}$ with minimum vanishing at $Q_{g-1}$ but is otherwise generic.
From Theorem \ref{modvbrc}, we obtain a stable vector bundle on the chain.

We construct in this way $k$ global sections of the vector bundle.
For $i=(k_1+1)j_1+j_2\le  (k_1+1)(g+k_1-d_1-1)$, the following sections  on  $C_i$ have order of vanishing on $P_i, Q_i$ adding up to $d_1$  
\[ s_{(j_2-1)r+1},\dots, s_{(j_2-1)r+r}, \  1\le j_2\le k_1;\ \ \ s_{(j_2-1)r+1},\dots, s_{(j_2-1)r+k_2},\  j_2= k_1+1\]

Recall that for the canonical limit linear series on a chain of elliptic curves, the line bundles on $C_i$ are  (Proposition \ref{can})
\[  L_i={\mathcal O}(2(i-1)P_i+2(g-i)Q_i) .\] 
  Therefore, the vector bundle on each of the elliptic components corresponding to $K\otimes E^*$ is well determined.
  There is then a limit linear series of rank $r$ degree $\bar d$ and dimension $\bar k$ with  
  \[ \bar d=r(2g-2)-d=r(2(g-1)+d_1-1)+r-d_2,\  \  \bar k=r(k_1-d_1+g-1)  \]
\begin{enumerate}[(a')]
 \item If  $j_2\not= k_1+1$,  the vector bundle on $C_i$ is  
 \[ {\mathcal O}_{C_i}((i+j_1-j_2)P_i+(2g+j_2-j_1-i-2)Q_i)^{\oplus r}.\]
  \item If  $j_2= k_1+1$, then ${\mathcal O}_{C_i}((i+j_1-j_2)P_i+(2g+j_2-j_1-i-2)Q_i)^{\oplus k_2}\oplus \bar L_1\oplus \dots\oplus \bar L_{r-k_2}$.

\item  For $  (k_1+1)(g+k_1-d_1-1)<i\le g-1$, the vector bundle on $C_i$, $ \bar L_1\oplus \dots\oplus \bar L_{r}$.

 \item On $C_g$ we get a generic  vector bundle of rank $r$ and degree $r(2g-3-d_1)+(r-d_2)$.
 \end{enumerate}

We construct in a similar  way $\bar k$ global sections of the vector bundle.
The following sections  on  $C_i$ have order of vanishing on $P_i, Q_i$ adding up to $d_1$, for $i=(k_1+1)j_1+j_2\le  (k_1+1)(g+k_1-d_1-1)$,   
\[ \bar s_{j_1r+1},\dots, \bar s_{j_1r+r}, \  1\le j_2\le k_1;\ \ \  \bar s_{j_1r+1},\dots, \bar s_{j_1r+k_2},\  j_2= k_1+1\]
 
Consider now the limit linear series corresponding to the bundle $E\otimes K\otimes E^*$.

 For  $i=(k_1+1)j_1+j_2,\ 0\le j_1\le g+k_1-d_1-2,\ 1\le j_2\le k_1$, the aspect of the limit linear series on $C_i$ has vector bundle  $({\mathcal O_{C_i}}(2(i-1)P_i+2(g-i)Q_i))^{r^2}$
 while for $ j_2= k_1+1$, the bundle is   $({\mathcal O_{C_i}}(2(i-1)P_i+2(g-i)Q_i))^{k_2^2}\oplus (\oplus_jL_j)$ with the $L_j$ line bundles of degree $2g-2$.
We have  sections $s_l\bar s_t, \ 1\le l\le k, 1\le t\le \bar k$ of the limit linear series  and our goal is to show that they are linearly independent.

Using Lemma \ref{serdistrdeg}, we will consider global sections of $E\otimes K\otimes E^*$ with degree $r$ on the first and last components and degree $2r$ on the intermediate ones.
The potentially non-zero sections on component $C_i$ must vanish at $P_i$  with order at least $2i-3, 2\le i\le g$ and at $Q_i$ with order at least $2g-2i-1, 1\le i\le g-1$ (see the proof of Lemma \ref{serdistrdeg}).

By construction, on the curve $C_i, i=(k_1+1)j_1+j_2,\ 0\le j_1\le g+k_1-d_1-2$, the following sections vanish to order $2i-2$ at $P_i$ and $2g-2i$ at $Q_i$ 
\[ 1\le j_2\le k_1, \ s_t\bar s _l, \  (j_2-1)r+1\le t\le j_2r  \ \ j_1r+1\le l\le (j_1r+1)r, \]
\[ j_2= k_1+1, \ s_t\bar s _l, \  k_1r+1\le t\le k_1r +k_2 \ \ j_1r+1\le l\le (j_1r+1)r, \]
Any other section vanishes to order at most two less at either $P_i$ or $Q_i$.
Consider then any possible linear dependence $\sum \lambda_{t,l}s_t\bar s _l=0$.
Restricting to the curve $C_i$ only the $s_t\bar s _l$ listed above are non-zero and these sections are sections of distinct line bundles in the direct sum. 
This forces $\lambda_{t,l}=0$, proving the independence of the $s_t\bar s _l$.
Therefore the image of the Petri map has dimension at least $k\bar k$.
Then, the Petri map is injective for curves and linear series close to the one we constructed.
\end{proof}


\section{Products of canonical and traceless endomorphisms }\label{secPoin}

Assume that we work in characteristic zero. Let $E$ be a vector bundle. Then,
\[ E^*\otimes E\cong {\mathcal O}\oplus Tr_0E\]
where $Tr_0E$ is the set of traceless endomorphisms.
The natural cup-product map
 \[ \psi:H^0(C,K)\otimes  H^0(C,K\otimes E^* \otimes E)\to H^0(C,K^2\otimes E^* \otimes E)\]
decomposes as direct sum of two maps
 \[  \psi_1: H^0(C,K)\otimes  H^0(C,K)\to H^0(C,K^2), \ \  \psi_2: H^0(C,K)\otimes  H^0(C,K\otimes Tr_0E)\to H^0(C,K^2\otimes Tr_0E)    \]
Then $\psi$ is onto if and only if both $ \psi_1,  \psi_2$ are onto.
The map $ \psi_1$ is onto if and only if $C$ is not hyperelliptic.
We prove here that $ \psi_2$ is onto for the general curve and general vector bundle.
It suffices to find a particular curve $X_0$ and a particular vector bundle $E_0$ on  $X_0$ with $h^0(X_0,K\otimes E_0^* \otimes E_0)=r^2(g-1)+1$ such that $\psi_2$ is onto for $X_0, E_0$.

\begin{definition}\label{defE0}
Let $X_0=C_1\cup\dots\cup C_g$ be a  generic chain of $g$ elliptic curves as in Definition \ref{chain}.
 Let $r,d$ be integers, with $1\le r, g\le d<g+r$.
  Let $E_0$ be a vector bundle on $X_0$ that restricts to an indecomposable vector bundle of rank $r$ and degree $1$ on $C_i, i=1,\dots, g-1$, 
while the restriction to $C_g$ is the direct sum of $h=gcd(r,d-g+1)$ generic indecomposable vector bundles of degree $\frac{d-g+1}h$ and rank  $\frac{r}h$.
Take all the gluing at the nodes to be generic.
  \end{definition}
  
  We show  that $h^0(C_0,K\otimes E_0^* \otimes E_0)=r^2(g-1)+1$ or equivalently, $h^0(C_0, E_0^* \otimes E_0)=1$:

  \begin{proposition}\label{hom}  Let $X_0,E_0$ be as in Definition \ref{defE0}.
  Then $Hom (E_0,E_0)=E_0^*\otimes E_0$    is a direct sum of line bundles of degree 0 on each component, only one of which  is trivial on $C_1,\dots, C_{g-1}$ while $h$ of them are trivial on $C_g$.
  The trivial line bundles  glue with each other in  $C_1,\dots, C_{g-1}$ and glue with one of the trivial line bundles on $C_g$ giving rise to a trivial line subbundle of $E_0^*\otimes E_0$ on $C_0$.
  The remaining gluing is generic.
  In particular, the limit linear series corresponding to $E_0^*\otimes E_0$ has a unique section.
  \end{proposition}
\begin{proof} We will use the additive and multiplicative structure of the moduli space of vector bundles on an elliptic curve (in characteristic zero) that was first  described in  \cite{A}.
Recall that if $gcd(r',d')=1$, two indecomposable vector bundles of rank $r'$ and degree $d'$ differ in tensorization with a line bundle of degree 0.
So, if  $F$ is such  an  indecomposable vector bundle, any other vector bundle of rank $r'$ and degree $d'$ is of the form $F\otimes L$ for some line bundle of degree 0.
Moreover, from the multiplicative structure for vector bundles on an elliptic curve (in characteristic zero) 
\[  F^*\otimes F=\oplus M_j, \  M_j \text{ line bundles of degree 0 of order dividing } r'.  \]
Therefore, on $C_i, i=1,\dots, g-1$, the vector bundles 
\[ E_i^*\otimes E_i=\oplus M_j^i, \  M_j^i \text{ line bundles of degree 0 on } C_i \text{ of order dividing } r,\ i=1,\dots, g-1.  \] 
 Let $F$ be an  indecomposable vector bundle of degree $d'=\frac{d-g+1}h$ and rank  $r'=\frac{r}h$ on $C_g$.
Then, there exist $h$ generic line bundles of degree 0 such that $E_g=F\otimes L_1\oplus \dots \oplus F\otimes  L_h$.
We can then compute $E_g^*\otimes E_g$:
\[ E_g^*\otimes E_g= (F^*\otimes L_1^*\oplus \dots \oplus F^*\otimes  L_h^*)\otimes  (F\otimes L_1\oplus \dots \oplus F\otimes  L_h)=F^*\otimes F\otimes (\oplus_{i,j}L_i^*L_j)=\]
\[ =\oplus M_j^g\otimes (\oplus_{i,j}L_i^*L_j), \  M_j^g \text{ line bundles of degree 0 on } C_g \text{ of order dividing } r',  \] 

Thinking of $E_0^*\otimes E_0$ as $Hom (E_0,E_0)$ , the identity morphism on $E_i, i=1,\dots, g-1$ glues with the identity morphism on $E_{i+1}$ . 
The identity morphism corresponds to the trivial subbundle of $ E_i^*\otimes E_i, i=1,\dots, g-1$ and to the diagonal in  ${\mathcal O}^h$ for $i=g$.
A non-trivial line bundle of degree 0 has no sections. 
Therefore, $E_0^*\otimes E_0$ has only one limit linear section.
\end{proof}

\begin{proposition}\label{onto}
With $E_0$ as in \ref{defE0}, there is a surjective product  map of limit linear series 
 \[ H^0(C_0,K)\otimes  H^0(C_0,K\otimes Tr_0E_0)\to H^0(C_0,K^2\otimes Tr_0E_0). \]
  \end{proposition}
  \begin{proof} From Proposition \ref{can}, the canonical series on $C_0$ has  line bundles $L_i={\mathcal O}(2(i-1)P_i+2(g-i)Q_i) $  on $C_i$
  and  limit sections $s_1,\dots, s_g$    with vanishing at $P_i, Q_i$ given respectively as 
        \[  \begin{matrix} i-2 & i-1 &\dots &2i-5&2i-4&2i-2&2i-1&\dots&g+i-2\\               2g-i-1&2g-i-2&\dots &2g-2i+2&2g-2i+1&2g-2i&2g-2i-2&\dots &g-i-1\\  \end{matrix}\]
We showed in our proof of \ref{hom}  that  the   bundle $Tr_0E_0$ has restrictions to $C_i,i=1,\dots, g-1$ of the form $\oplus_{k=1\dots r^2-1}L_k$ where the $L_k$ are line bundles of degree 0.
Moreover, by the genericity of the chain of elliptic curves, $P_i-Q_i$ is not a torsion point in the group structure of the elliptic curve,
Therefore,  none of the line bundles appearing in the decomposition of  $Tr_0E_i$ is of the  form ${\mathcal O}(aP_i+bQ_i), a+b=0$ for any $a,b$.
On  $C_g $, we still have a direct sum of line bundles of degree 0 but there $h-1$ of them are the trivial bundle while the rest are not  of the  form ${\mathcal O}(aP_i+bQ_i), a+b=0$ for any $a,b$.

Therefore, the   bundle $K\otimes Tr_0E_0$ has restrictions to $C_i,i=1,\dots, g-1$ of the form $\oplus_{k=1\dots r^2-1}\bar L_k$.
 Here the $\bar L_k$ are line bundles of degree $2g-2$  not of the form ${\mathcal O}(aP_i+bQ_i), a+b=2g-2$.
 On $C_g$, $h-1$ of the line bundles are  ${\mathcal O}((2g-2)Q_i)$ while the rest are generic.

This allows us to find $r^2-1$ sections  $ t^i_{1,j},\dots,t^i_{r^2-1,j}$ of the restriction of  $Tr_0E_0$  to $C_i$ vanishing with order   $i+j-2$ at $P_i$ and $2g-i-j-1$ at $Q_i$ $j=1,\dots,g-1$
spanning the fibers of the bundle at both $P_i,Q_i$ (see Lemma \ref{lemseclselc}).
On $C_g$, for $j=g-1$, $h-1$ of the sections $ t^g_{1,g-1},\dots,t^g_{r^2-1,g-1}$ have vanishing order $2g-2$ instead of the minimal  required of $2g-3$ at $P_g$.

As  $ t^i_{1,j},\dots,t^i_{r^2-1,j}$ span  the fibers of $K\otimes Tr_0E_0$ at $P_i,Q_i$,
we can  assume that $t_{k,j}^i$ glues  with  $t_{k,j}^{i+1}$  in the identification of $Q_i$ and  $P_{i+1}$.
 
  We will show that in the cup product  $H^0(C_0,K)\otimes  H^0(C_0,K\otimes Tr_0E_0)\to H^0(C_0,K^2\otimes Tr_0E_0 )$, the images of 
   \begin{equation}\label{colsec} s_1\otimes t_{k,1}, s_1\otimes t_{k,2},s_1\otimes t_{k,3},s_2\otimes t_{k,2}, s_2\otimes t_{k,3},s_2\otimes t_{k,4},\dots,
   s_{g-3}\otimes t_{k,{g-3}}, s_{g-3}\otimes t_{k,{g-2}},s_{g-3}\otimes t_{k,{g-1}}, \end{equation}
 \[   s_{g-2}\otimes t_{k,{g-2}},s_{g-2}\otimes t_{k,{g-1}} ,  s_{g-1}\otimes t_{k,{g-2}},s_{g-1}\otimes t_{k,{g-1}}  , s_{g}\otimes t_{k,{g-2}},s_{g}\otimes t_{k,{g-1}}  ,\ k=1,\dots r^2-1\]
    are linearly independent.
    
Using Lemma \ref{serdistrdeg}, we will consider the distribution of degrees that assigns degree  $3(r^2-1) $ to the components $C_1, C_{g-2}, C_{g-1}, C_{g} $ and degree $4(r^2-1) $ to the remaining components.
That is, we consider the vector bundles 
\[ E_1(-(4g-7)Q_1);\  E_i(-(4i-5)P_i-(4g-4i-3)Q_i),  i=2,\dots, g-3, \]
\[ E_{g-2}(-(4g-13)P_{g-2}-6Q_{g-2}),\ E_{g-1}(-(4g-10)P_{g-1}-3Q_{g-1}), \ E_{g}(-(4g-7)P_{g})  \]
 Assume that we had a linear combination of the chosen sections that equals 0.

 On $C_1$, the only sections in (\ref{colsec})  that vanish to  order at least  $4g-7$  at $Q_1$ are 
  \[ s_1\otimes t_{k,1}, s_1\otimes t_{k,2},s_1\otimes t_{k,3}, k=1,\dots, r^2-1\]
    vanishing to orders  $4g-5, 4g-6, 4g-7$ respectively.
    The remaining sections, as sections from the series with adjusted degrees, are identically zero on $C_1$.
    Hence, the initial linear dependence condition among all the sections in (\ref{colsec}) reduces on $C_1$ to a linear dependence condition among the 
   $s_1\otimes t_{k,1}, s_1\otimes t_{k,2},s_1\otimes t_{k,3}$ only.
   The orders of vanishing at $P_1$ for these three sets of sections are $0,1,2$ respectively. 
   The images  of $s_1\otimes t_{k,j}$ for fixed $j$ and varying $k$ are linearly independent as they map onto the fiber of $K\otimes E^*\otimes E$ at $Q_1$.
     Given that the three sets have different orders of vanishing at $P_1$ and that each set is linearly independent,
     $s_1\otimes t_{k,1}, s_1\otimes t_{k,2},s_1\otimes t_{k,3}$ have $0$ coefficients on any zero  linear combination of  list (\ref{colsec}).
    
     Let us assume by induction that for the sections  $s_1\otimes t_{k,1}, s_1\otimes t_{k,2},s_1\otimes t_{k,3}\dots  s_{l-1}\otimes t_{k,l-1}, s_{1-1}\otimes t_{k,l},s_{l-1}\otimes t_{k,l+1}$ 
 the coefficients   in the trivial linear combination are zero.
  Look now at the restriction  to $C_{l}$ of  the sections in (\ref{colsec}).
  Among the  sections with potentially non-zero coefficients,  the only ones that vanish at $Q_l$ to  order at least  $4g-4l-3$  are 
  $s_l\otimes t_{k,l}, s_l\otimes t_{k,l+1},s_l\otimes t_{k,l+2}$ that vanish to orders  $4g-4l-1, 4g-4l-2, 4g-4l-3$ respectively.
  The remaining sections, as sections of the vector bundle with adjusted degrees vanish identically on $C_l$.
  The orders of vanishing at $P_l$ for each of these three sets of sections are different while the images on a stalk for a fixed $j$ but varying the $k$ are linearly independent.
  Therefore, the coefficients on the linear combination are zero.
  
   A similar reasoning works on the last three components where only $2(r^2-1)$ of the sections in (\ref {colsec}) are not identically zero 
  and those have different orders of vanishing and map onto the fibers.
    \end{proof}
    
 \section{Funding and/or Conflicts of interests/Competing interests.}
 
The authors have no relevant financial or non-financial interests to disclose.
 The author received partial support from NSF  Grant Account Number: 104301-0000


\end{document}